\documentclass[letterpaper, 10 pt, conference]{ieeetran}  
\IEEEoverridecommandlockouts                  
\usepackage{graphicx}  
\usepackage{subfigure}
\usepackage{multirow}
\usepackage{textcomp}
\usepackage{amsmath}
\usepackage{ifthen}
\usepackage{amssymb}
\usepackage{amsfonts}
\usepackage{subcaption}
\usepackage{colortbl}
\usepackage{algorithm,comment}
\usepackage{algpseudocode}
\usepackage{hyperref}
\usepackage{mathtools,color}
\usepackage{xcolor}
\usepackage{array}
\newtheorem{theorem}{\textbf{Theorem}}
\newtheorem{remark}{\textbf{Remark}}
\newtheorem{corollary}{\textbf{Corollary}}
\newtheorem{lemma}{\textbf{Lemma}}
\newtheorem{definition}{\textbf{Definition}}
\newtheorem{assumption}{\textbf{Assumption}}

\title{\LARGE \bf
A Control Barrier Function Approach to Constrained Resource Allocation Problems in a Maximum Entropy Principle Framework
}

\author{Alisina Bayati$^{\dagger}$, Dhananjay Tiwari$^{\dagger}$, and Srinivasa Salapaka
\thanks{$^{\dagger}$These authors contributed equally to this work.}
\thanks{All three authors are from the University of Illinois, Urbana-Champaign. Emails: \texttt{\{abayati2, dtiwari2, salapaka\}@illinois.edu}. We acknowledge the support of National Aeronautics and
Space Administration under Grant NASA 80NSSC22M0070 for this work.} }

\newcommand{\bb}[1]{\mathbb{#1}}

\newcommand{\mcal}[1]{\mathcal{#1}}
\newcommand{\rcb}[1]{\left[#1\right]}
\newcommand{\cb}[1]{\left\{#1\right\}}
\newcommand{\rb}[1]{\left(#1\right)}
\newcommand{\abs}[1]{\left|#1\right|}
\newcommand{\mb}[1]{\mathbb{#1}}

\begin{document}

\maketitle
\thispagestyle{empty}
\pagestyle{empty}

\begin{abstract}
This paper presents a novel approach to solve capacitated facility location problems (FLP) that encompass various resource allocation problems. FLPs are a class of NP-hard combinatorial optimization problems, involving optimal placement and assignment of a small number of facilities over a large number of demand points, with each facility subject to upper and lower bounds on its resource utilization (e.g., the number of demand points it can serve). To address the challenges posed by inequality constraints and the combinatorial nature of the solution space, we reformulate the problem as a dynamic control design problem, enabling structured constraint handling and enhanced solution efficiency. Our method integrates a Control Barrier Function (CBF) and Control Lyapunov Function (CLF)-based framework with a maximum-entropy principle-based framework to ensure feasibility, optimality, and improved exploration of solutions. Numerical experiments demonstrate that this approach significantly enhances computational efficiency, yielding better solutions and showing negligible growth in computation time with problem size as compared to existing solvers. These results highlight the potential of control-theoretic and entropy-based methods for large-scale facility location problems.
\end{abstract}

\vspace{1em}

\section{Introduction}
Constrained resource allocation problems span various domains, including bandwidth allocation in cognitive radio networks \cite{awoyemi2016solving}, cloud resource provisioning \cite{cloud_comp_vinothina2012survey}, LLM task allocation in cloud-edge networks \cite{llm_he2024large}, 
data routing in 5G networks \cite{small_cell_wu2017qoe}, air-pathway clustering for AAM \cite{li2022traffic}, UAV scheduling \cite{tiwari2025multi}, and supply-chain logistics \cite{supply_chain_griffis2012metaheuristics}. Despite their varied goals, these problems can be framed as variants of the {\em Facility Location Problem (FLP)}—which involves determining the optimal placement of $M$ facilities (e.g., warehouses, service centers) to serve $N$ demand points while minimizing costs such as transportation and operations. The key decision variables are {\em assignments}, specifying which demand point is allocated to which facility, and {\em facility locations}, representing the spatial or resource attributes of facilities.  The differences in these problems arise from application specific distance functions and constraints, such as AAM scheduling, which prioritizes shortest routes while avoiding collisions, and cloud computing, which focuses on fast computations with equitable task distribution across resources.

These mixed-integer problems, are combinatorial with the number of decision variables growing exponentially with problem size. For example, a simple unconstrained FLP with $M$ facilities and  $N\gg M$ demand points results in $2^{NM}$ binary assignment variables and $M$ variables in ${\bb{R}}^2$. Such problems are NP-hard, with cost surfaces containing numerous poor local minima \cite{gray1982multiple}.  Various approaches address these problems, including classical optimization \cite{classical_opt_CRN_wang2010robust}, integer constraint relaxation \cite{integer_constraint_li2013efficient}, and heuristics/meta-heuristics such as $k$-means clustering \cite{original_kmeans_lloyd1982least}, genetic algorithms \cite{clustering_GA_cole1998clustering}, simulated annealing \cite{SA_book_van1987simulated}, and tabu search \cite{tabu_glover2007principles}. However, these methods are often problem-specific and prone to suboptimal solutions that are highly sensitive to initialization \cite{rose1991deterministic}.

We highlight the Deterministic Annealing (DA) algorithm \cite{rose1991deterministic}, an optimization technique based on the Maximum Entropy Principle (MEP). 
DA has demonstrated superior efficiency over many heuristic and meta-heuristic methods, often achieving faster computation and lower objective costs \cite{ourAIAA, srivastava2020simultaneous}. 
It has been successfully applied to clustering, combinatorial optimization, and machine learning, serving as a robust alternative to traditional heuristics. 
Notable applications include
airship deployment in wireless communication \cite{airships_hap_wang2011energy}, training neural-networks \cite{learning_mavridis2022online, supervised_learning_gupta2019noisy}, multi-robot path planning \cite{DA_APF_dai2022robot} and task allocation \cite{robot_task_alloc_david2022deterministic}, 
traveling salesman problem \cite{baranwal2017multiple}, MDPs and reinforcement learning \cite{srivastava2021parameterized}.

In DA approach, combinatorial binary assignment variables are replaced by soft probability distributions, and the cost objective is reformulated as a {\em free-energy} function, combining the relaxed cost objective with Shannon entropy to quantify uncertainty. DA solves a sequence of optimization problems parameterized by an {\em annealing parameter} $\beta$, which controls the balance between entropy and the original cost function. When $\beta$ is small, the free energy function is convex, allowing its global minimum to be efficiently located.
As $\beta$ increases, free energy is minimized iteratively, initializing each step with the previous solution. 
As $\beta \to \infty$, the relaxed problem takes the form of original NP-hard problem. This controlled exploration enables DA to efficiently balance exploration and exploitation, guiding the solution toward global or near-global optima.

While the MEP-based framework is efficient in cost and computation time, its application to scalable constrained resource allocation, particularly with capacity constraints limiting facility assignments, remains underexplored.
Works such as \cite{baranwal2017clustering, srivastava2022inequality} address constrained FLP variants like FLPO and LMDP. \cite{baranwal2017clustering} enforces capacity constraints as equalities within DA, solving for auxiliary capacity variables but overlooking inequality constraints. 
In contrast, \cite{srivastava2022inequality} handles capacity constraints by formulating them as inequality constraints and augmenting the free energy with a penalty term for each constraint. 
However, the effectiveness of this approach is highly sensitive to the choice of penalty functions and their relative weighting with respect to the free energy, making it prone to constraint violations or suboptimal solutions. 
Alternatively, solvers such as Sequential Least Squares Programming (SLSQP) and \textit{trust-constr} can handle inequality constraints in large-scale nonlinear optimization. While both methods solve quadratic subproblems, SLSQP~\cite{boggs1995sequential} uses second-order Lagrangian approximations with linearized constraints and \textit{trust-constr}~\cite{byrd1999interior} uses an interior-point method with trust-regions. 
However, in our applications~\cite{tiwari2025multi,ourAIAA}, these methods exhibit poor runtime performance, primarily due to the computational overhead associated with second-order Hessian evaluations.

This article presents a systematic and efficient method for incorporating equality and inequality constraints in FLPs. The approach builds on Control Barrier Functions (CBFs) and Control Lyapunov Functions (CLFs), originally developed for nonlinear control-affine systems in applications such as robotics and UAVs, where safety and stability are critical. The framework in~\cite{8796030, 6760327} ensures constraint satisfaction while steering the system toward an equilibrium defined by the CLF. Under mild conditions, the resulting controller is Lipschitz continuous~\cite{6760327, 7782377}.

We recast our {\em static} constrained FLP as a {\em dynamic} control design problem, and apply the CLF-CBF based framework to  solve it. First, we relax the problem by replacing binary assignments with smooth probability distributions and defining the Free Energy function ${\mathcal{F}}$ parameterized by the annealing parameter $\beta$, following the DA approach. 
At each $\beta$, we {\em introduce} control-affine dynamics $\dot z_{\beta}= q_0(z_{\beta})+q_1(z_{\beta})u, z_{\beta}(0)=z_{\beta,0}$, where the state $z$ comprises the decision variables (assignment distributions and facility locations), and $u$ represents the control input. 
The objective is to steer the system toward a local minimum while ensuring feasibility. This is achieved by designing $u$ via a quadratic program (QP) that enforces stability (via a CLF) and safety (via CBFs). We establish that (i) $\mcal F$ acts as a CLF, ensuring the existence of a feasible control $u$ such that $\dot {\mathcal{F}}\leq 0$, and (ii) the closed-loop system maintains feasibility for well-posed problems and converges to stationary points satisfying the Karush-Kuhn-Tucker (KKT) conditions of the relaxed constrained FLP. Annealing is incorporated by iterating this process for increasing $\beta_k$ and reinitializing $z_{\beta_{k^+},0}$ with the stationary point from $\beta_{k-1}$.

We compared our implementation against three baselines: (1) Safe Gradient Flow, a CBF-based method for nonlinear programming \cite{sgf_allibhoy2023control}, (2) SLSQP from \texttt{SciPy} \cite{virtanen2020scipy} and (3) DA penalty-based approach from \cite{srivastava2022inequality}. Our approach achieves comparable solution quality while being nearly 20 times and 240 times faster as compared to SGF and SLSQP, respectively. This is primarily due to efficient QP solved compared to Safe Gradient Flow and the elimination of costly Hessian computations required by SLSQP \cite{boggs1995sequential}. 
Further, it shows better constraint handling as compared to DA-based penalty approach, due to a structured and problem independent solution approach.
A key advantage of our framework is its adaptability to dynamic settings. With minor modifications to the QP governing the control input, it can handle time-varying optimization landscapes where demand points evolve according to known dynamics. This allows our control design to continuously track the local minimum of the evolving free energy surface. We are exploring this direction for applications such as battleship-based surveillance and multi-robot dynamic coverage, where performance optimization and safety constraints, such as collision avoidance, must be jointly addressed.

\section{Problem formulation} \label{sec:RA_formulation}
Consider an FLP in the context of supply-chain logistics, where \( N \) demand points (e.g., commodity production centers) are located at \( x_i \in \mathbb{R}^d \) for \( 1 \leq i \leq N \). The goal is to assign \( M \ll N \) facilities—hereafter referred to interchangeably as resources—located at \( y_j \in \mathbb{R}^d \) for \( 1 \leq j \leq M \), to serve these demand points so as to minimize the average transportation cost. This leads to the following optimization problem:

\(\textbf{P}_{\textbf{flp}}\)\textbf{:}
\vspace{-15pt}
\begin{subequations} \label{prob:P_flp}
\begin{align}
    \min_{\substack{y_j \in \mathbb{R}^d \\ \nu_{j|i} \in \{0,1\}}} 
    &\mathcal{D} :=  \sum_{i=1}^{N} p_i \sum_{j=1}^{M} \nu_{j|i} \, d(x_i, y_j) \label{def:D_Total_Cost} \\
    \text{s.t.} ~~~
    & \sum_{j=1}^M \nu_{j|i} = 1, \quad \forall\, i = 1, \dots, N, \label{probcon} \\
     L_j \leq &\sum_{i=1}^{N} p_i \, \nu_{j|i} \, c_{ij} \leq C_j, \quad \forall\, j = 1, \dots, M. \label{capcon}
\end{align}
\end{subequations}
Here $p_i>0$ represents a priori known relative weight or importance of the $i$th demand point, with $\sum_{i=1}^N p_i=1$. The  binary variable $\nu_{j|i}$ indicates assignment ($\nu_{j|i}=1$ when resource $j$ serves demand point $i$, otherwise $0$). The function $ d: \mathbb{R}^d \times \mathbb{R}^d \to \mathbb{R}_{\geq 0},$ represents   the  cost associated with assigning $x_i$ to $ y_j$ (e.g., squared Euclidean distance $ d(x_i, y_j) = \|x_i - y_j\|_2^2$); therefore $\mathcal{D}$ represents the average cost over all assignments. Constraint (\ref{probcon}) ensures that each demand point $i$ is serviced by exactly one resource.

Constraint~\eqref{capcon} enforces the capacity limits at each facility, where \( c_{ij} \) denotes the resource consumption of the \( i \)th demand point when assigned to the \( j \)th facility. The bounds  $L_j$ and $C_j$ represent the minimum and maximum resource usage allowed for $j$th facility respectively. The interpretation of $c_{ij}$ varies by  the application: In \textit{manufacturing}, it represents the machine workload or processing effort required to complete task $i$ at workstation $j$ 
\cite{kovacs2021utilizing};  In \textit{cloud computing}, it quantifies the computational resources—such as CPU, memory, or I/O bandwidth needed to allocate task $i$ to server $j$ 
\cite{perez2022dynamic}.
The upper bound $C_j$ prevents overloading or overcommitment, while the lower bound $L_j$ ensures efficient resource use, promoting load balancing—crucial for applications like cloud computing and power distribution.

Without the capacity constraints \eqref{capcon}, \(\text{P}_{\text{flp}}\) reduces to the {\em unconstrained FLP}, which is solvable via MEP. Since our approach builds on the resulting DA algorithm, we first introduce it for the unconstrained case.

\section{Maximum Entropy Principle (MEP) based solution approach}

\subsection{{Unconstrained Facility Location Problem}}\label{subsec:MEP_uncons_RA}

The MEP-based DA algorithm solves a related parameterized optimization problem, defined as below:

\(\textbf{P}_{\textbf{unconstr}}(\boldsymbol{\beta})\)\textbf{:}
\vspace{-10pt}
\begin{subequations}\label{prob:P_unconstr_beta}
\begin{align}
    \min_{{y_j^\beta \in {\bb{R}}^d},{p^\beta_{j|i}\in [0,1]}}
    & \mcal{F^\beta } :=\mathcal{D^\beta}-\frac{1}{\beta}\mathcal{H^\beta}\label{cfpunconst}\\
    \text{s.t.} \quad \quad  & \sum_{j=1}^M p^\beta_{j|i}=1, \quad \forall 1 \leq i \leq N,
\end{align}
\end{subequations}
where $\beta$ is an {annealing} parameter. The algorithm leverages the algebraic structure of solutions of $\textrm{P}_{\textrm{unconstr}}\rb \beta$ to obtain an approximate local minimum of the unconstrained FLP. 

Here $\mcal{D^\beta} := \sum_{i} p_i \sum_j p^\beta_{j|i} d(x_i, y_j),$ represents a relaxation of cost function $\mcal{D}$ in (\ref{def:D_Total_Cost}) of the unconstrained FLP, where we have replaced 
the binary association variables $\nu_{j|i}$ with a set of soft probability associations $p^\beta_{j|i} \in \rcb{0,1}.$
Further, $\mathcal{H^\beta} :=- \sum_{i} p_i \sum_j p_{j|i}^\beta \log p_{j|i}^\beta$ is the (conditional) Shannon entropy of the distribution $\{p_{j|i}^\beta\}$, that measures its uncertainty, 
with higher entropy indicating a more evenly distributed association of demand point $i$ across all $M$ resources. Such associations typically result in solutions to the
$\textrm{P}_{\textrm{unconstr}}\rb \beta$
that are less sensitive to initializations. However such associations imply higher values of $\mcal{D^\beta}$ since $i$th demand point is not associated to the nearest (local) resource location.  The free-energy function \( \mathcal{F}^\beta \) in~\eqref{cfpunconst} thus captures a trade-off between approximation accuracy and sensitivity to initialization.
At low $\beta$, $\{p^\beta_{j|i}\}$ has high entropy, reducing sensitivity to initialization, while at high $\beta$, entropy decreases, leading to more localized associations, and better approximations to the unconstrained FLP. 

Applying the first-order necessary conditions to the unconstrained Lagrangian associated with \( \textrm{P}_{\textrm{unconstr}}(\beta) \) yields:
\begin{equation} \label{eq:unconstrained_updates}
    p_{j|i}^\beta = \frac{e^{-\beta d(x_i, y_j^\beta)}}{\sum_{\ell} e^{-\beta d(x_i, y_{\ell}^\beta)}} ~~ \forall i,j, \quad
    y_j^\beta = \frac{\sum_i p_i p^\beta_{j|i} x_i}{\sum_i p_i p^\beta _{j|i}} ~~ \forall j,
\end{equation} 
where $p^\beta_{j|i}$ corresponds to the softmax (Gibbs) distribution for a given set of $\{y_j^\beta\}$ and the expression of $y_j^\beta$ represents the weighted centroid of the assigned demand points for given $\{p^\beta_{j|i}\}.$ Thus the solution at any $\beta$ is obtained by solving the implicit equations \eqref{eq:unconstrained_updates}, which is shown to converge to a fixed-point \cite{rose1991deterministic, baranwal2017clustering}. 
Note that at $\beta=0$, the solutions satisfy $p^\beta_{j|i}=1/M$ $\forall i, j$, and $y^\beta_j=\sum_i p_ix_i$, which is independent of the initialization and a unique global solution to $\textrm{P}_{\textrm{unconstr}}\rb \beta$. While at $\beta=\infty$, the solutions satisfy, $p^\beta_{j|i}=1$, if resource $j$ is a closest to demand point $i$, else it is $0$. 
Thus the corresponding solution is a  (local) solution to the unconstrained FLP. 
Also at $\beta=\infty$, the $\mcal{F^\beta}=\mcal{D^\beta}$ since $\frac{1}{\beta}H^\beta=0$ and $\mcal{D^\beta}$ has the same form as $\mcal{D}$. 

{\em \bf Annealing Process:} The DA algorithm implements an {\em annealing} schedule by iteratively solving the optimization problem \( \textrm{P}_{\textrm{unconstr}}(\beta_k) \) over a sequence of increasing parameters \( \{ \beta_k = \gamma^k \beta_0 \} \), where \( \gamma > 1 \) and \( \beta_0 > 0 \) is small. At each step, the solution from \( \beta_{k-1} \) is used to initialize the problem at \( \beta_k \). In the early stages, when \( \beta \approx 0 \), the resource locations \( \{ y_j^\beta \} \) tend to be similar, typically concentrated near the weighted (by \( p_i \)) centroid of the demand points. As \( \beta \) increases, the resources progressively separate into {\em distinct} clusters aligned with the structure of the demand. The probability assignments \( p^\beta_{j|i} \) eventually converge to hard associations \( \nu_{j|i} \). Thus, DA can be viewed as identifying a global minimum at low \( \beta \), where \( \mathcal{F} \) is convex, and then tracking this minimum as \( \beta \) increases. 

\subsection{Constrained Facility Location Problem}
Here we modify the MEP based approach to include the capacity constraints. Accordingly we formulate a parameterized $\textrm{P}_{\textrm{constr}}(\beta)$ and develop a similar annealing process, where we solve repeatedly $\textrm{P}_{\textrm{constr}}(\beta_k)$ at a sequence of parameters $\{\beta_k=\gamma^k\beta_0\},$  using the solution from $\beta_{k-1}$ to initialize  $\textrm{P}_{\textrm{constr}}(\beta_k)$. Here $\textrm{P}_{\textrm{constr}}(\beta)$ is given by:

\(\textbf{P}_{\textbf{constr}}(\boldsymbol{\beta})\)\textbf{:}
\vspace{-10pt}
\begin{subequations}\label{prob:P_constr_beta}
    \begin{align}
        \min_{{y_j^\beta \in {\bb{R}}^d, p^\beta_{j|i}\in [0,1]}}
        & \mcal{F^\beta }:=\mathcal{D^\beta}-\frac{1}{\beta}\mathcal{H^\beta}\label{cfpconst}\\
        \text{s.t.} \quad \quad & L_j \leq \sum_i p_i \, p^\beta_{j|i} \, c_{ij} \leq C_j, \quad \forall j,\label{relconstr}\\
        \ \ & \sum_{j=1}^M p^\beta_{j|i}=1, \quad \forall i,
    \end{align}    
\end{subequations}
where we have relaxed the capacity constraints by using probability distributions $p^\beta_{j|i}$ instead of hard binary assignments $\nu_{j|i}$. The main issue that arises here is that DA algorithm cannot easily be extended to accommodate inequality constraints (\ref{relconstr}). Such attempts, as described in the introduction section, are either computationally inefficient or difficult to implement in practice.   

\begin{remark}
Unlike \( \textrm{P}_{\textrm{unconstr}}(\beta) \), the probability associations \( \{p_{j|i}^\beta\} \) in \( \textrm{P}_{\textrm{constr}}(\beta) \) may not converge to binary assignments for all demand points as \( \beta \to \infty \). This is because the nearest resource may lack sufficient capacity to fully serve a given demand point, resulting in fractional allocations even in the zero-entropy limit. In such cases, a probabilistic interpretation is useful: the values \( \{p_{j|i}^\beta\} \) represent allocation probabilities across multiple resources. From this perspective, the expected resource consumption satisfies the given lower and upper bounds, even if individual realizations may slightly violate them.
\end{remark}

In the following section, we introduce a CBF-CLF-based algorithm designed to generate trajectories that asymptotically converge to the KKT points of general nonlinear optimization problems involving both equality and inequality constraints. We will later show that \( \textrm{P}_{\textrm{constr}}(\beta) \) arises as a special case within this framework. 

\section{A Control-Theoretic Viewpoint on Optimization}

In this section, we explore how a general nonlinear optimization problem can be reformulated as a control design problem, where the goal is to steer the decision variable from a feasible initial condition to a stationary point satisfying the Karush-Kuhn-Tucker (KKT) conditions. 
We propose an approach in which control inputs directly act on the decision variable, ensuring both descent toward stationarity and constraint satisfaction along the trajectory.

As an alternative viewpoint, we reference the Safe Gradient Flow (SGF) approach introduced in~\cite{sgf_allibhoy2023control}, where control inputs influence the Lagrange multipliers associated with the constraints of the problem. In this formulation, the system flows along the gradient of the objective function, and control inputs intervene when feasibility is near-violation.

We detail our approach in the next subsection and briefly summarize the SGF method, referring the reader to \cite{sgf_allibhoy2023control} for further information. 

\subsection{Our Approach}

Consider the following optimization problem, where the minimum exists:
\begin{align}
    &\min_{z \in \mathbb{R}^b} \quad f(z) \ \ \text{s.t.} \ \ \begin{cases}
     g_i(z) = 0, \quad i = 1, 2, \dots, n,\\[1mm]
    h_j(z) \geq 0, \quad j = 1, 2, \dots, k.
    \end{cases}
    \label{eq:opt_primal}
\end{align}
As discussed, we reinterpret \eqref{eq:opt_primal} as a control design problem governed by the dynamics
\begin{align} \label{eq:z_dynamics}
    \dot{z}(t) = u(t), \quad z(0) = z_0,
\end{align}
where \(u(t)\) is designed to steer the state \( z(t) \) from a feasible initial point \(z_0\) toward a stationary (KKT) point. While KKT points do not necessarily correspond to local or global minima, our primary objective is to ensure convergence to stationarity. Further refinement—such as escaping saddle points—can be addressed using perturbation-based techniques after convergence is attained.

Before presenting Theorem~\ref{Thrm:main_theorem}, the main result of this work, we first introduce the following definitions and assumptions.

\begin{definition}
Define the sets \(\mathcal{C}\) and \(\mathcal{S}\) as the regions satisfying the equality and inequality constraints, respectively:
\begin{align} 
    \mathcal{C} &= \{z \in \mathbb{R}^b \mid g_i(z) = 0, \quad \forall i = 1, \dots, n\}, \label{eq:feasible_sets_equality} \\  
    \mathcal{S} &= \{z \in \mathbb{R}^b \mid h_j(z) \geq 0, \quad \forall j = 1, \dots, k\}. \label{eq:feasible_sets_inequality}
\end{align}
\end{definition}

\begin{definition}
At any point \( z \in \mathcal{C} \cap \mathcal{S} \), a vector \( \eta \in \mathbb{R}^b \) is called a \underline{feasible direction} if it satisfies:
\begin{align}
\begin{cases}
\left\langle \nabla g_i(z), \eta \right\rangle = 0, &\quad \forall i = 1, \dots, n,\\
\left\langle \nabla h_j(z), \eta \right\rangle \geq 0, &\quad \forall j \text{ such that } h_j(z) = 0,
\end{cases}
\end{align}
where \( \langle \cdot, \cdot \rangle \) denotes the standard inner product in \( \mathbb{R}^b \).
\end{definition}

Assume that the feasible set \( \mathcal{C} \cap \mathcal{S} \) is nonempty and path-connected. Without loss of generality, suppose \( f(z) \geq 0 \) for all \( z \in \mathcal{C} \cap \mathcal{S} \); otherwise, a constant shift can be applied to ensure nonnegativity. Additionally, suppose the following assumptions hold:

\begin{assumption} \label{assump:constraint_gradients}
The functions \( f \), \( g_i \) for all \( i = 1, \dots, n \), and \( h_j \) for all \( j = 1, \dots, k \) are $C^1$ with locally Lipschitz gradients. Furthermore, every point \( z \in \mathcal{C} \cap \mathcal{S} \) is \underline{regular}, meaning that the gradients of all equality constraints and active inequality constraints are linearly independent at $z$.
\end{assumption}

\begin{assumption} \label{assump:coercivity}
Either the set \( \mathcal{S} \) is bounded, or the objective function \( f(z) \) is \underline{coercive} on \( \mathcal{C} \, \cap \, \mathcal{S} \), meaning:
\[
\lim_{{z \in \mathcal{C}\cap\mathcal{S}, \ \|z\|\to\infty}} f(z) = \infty.
\]
\end{assumption}

\begin{theorem} \label{Thrm:main_theorem}
Under the assumptions \ref{assump:constraint_gradients} and \ref{assump:coercivity}, consider the feedback control \( u^*(z) \) defined as the solution to the following quadratic program (QP):
\begin{subequations} \label{eq:QP_general}
\begin{align} 
    u^*(z) &:= \arg\min_{(u,\delta) \in \mathbb{R}^{b+1}} \|u\|^2 + q \, \delta^2 \label{eq:QP_general_obj} \\
    \text{s.t.}  ~~~
    &\dot{f}(z;u) := \langle \nabla f(z), u \rangle \leq -\gamma(f(z)) + \delta, \label{eq:QP_general_fdot} \\
    &\dot{h}_j(z;u) := \langle \nabla h_j(z), u \rangle \geq -\alpha_j(h_j(z)), 
    &&\forall j, \label{eq:QP_general_hdot} \\
    &\dot{g}_i(z;u) := \langle \nabla g_i(z), u \rangle = 0, 
    &&\forall i, \label{eq:QP_general_gdot}
\end{align}
\end{subequations} 
where \(\gamma\) is a class~\(\mathcal{K}\) function, each \(\alpha_j\) is an extended class~\(\mathcal{K}_\infty\) function and \(q > 0\) is a  constant. Then, the following properties hold:
\begin{enumerate} 
    \item For any \( z \in \mathcal{C} \cap \mathcal{S} \), the control \( u^*_z : = u^*(z) \) defined by \eqref{eq:QP_general} exists uniquely and ensures \( \dot{f}(z; u_z^*) \leq 0 \), with strict inequality \( \dot{f}(z;u^*_z) < 0 \) if and only if \( z \) is not a KKT point, and equality \( \dot{f}(z;u^*_z) = 0 \) if and only if \( z \) satisfies the KKT conditions.
    \item \(u^*(z)\) is locally Lipschitz continuous on \(\mathcal{C} \cap \mathcal{S}\).
    \item The trajectory \( z(t) \) remains in  \( \mathcal{C} \cap \mathcal{S} \) for all \( t \geq 0 \), and asymptotically converges to a KKT point of \eqref{eq:opt_primal}.
\end{enumerate}
\end{theorem}

The core idea of Theorem~\ref{Thrm:main_theorem} is to interpret \( f \) (or any smooth, monotonic transformation such as \( f^2 \)) as a CLF-like function that drives the system toward stationarity, while the inequality constraints \( h_j \) are treated as CBFs—following the framework of Ames et al.~\cite{8796030, 7782377}—to ensure forward invariance of the set \( \mathcal{S} \).

\begin{proof}[Theorem~\ref{Thrm:main_theorem}]

\subsubsection*{Part 1}

To see this, note that \eqref{eq:QP_general} defines a strictly convex QP, and the feasible set over \( (u, \delta) \in \mathbb{R}^{b+1} \) is nonempty, as the trivial solution \( u = 0 \), \( \delta = \gamma(f(z)) \) satisfies all constraints for any \( z \in \mathcal{C} \cap \mathcal{S} \). Therefore, \( u^*(z) \) exists and is unique for all such \( z \). Moreover, any nonzero control $u$ satisfying \eqref{eq:QP_general_hdot} and \eqref{eq:QP_general_gdot} that leads to \( \dot{f}(z; u) > 0 \) incurs a strictly higher cost than the trivial solution, which yields \( \dot{f}(z; u) = 0 \). Hence, \( \forall z \in \mathcal{C} \cap \mathcal{S} \), $u_z^*$ renders \( \dot{f}(z;u^*_z) \leq 0 \).

If \( z \in \mathcal{C} \, \cap \, \mathcal{S} \) is not a KKT point, we claim (and later prove) that there exists a control \( \tilde{u} \) satisfying \eqref{eq:QP_general_hdot} and \eqref{eq:QP_general_gdot} and yields \( \dot{f}(z;\tilde u) < 0 \). It then follows that, for sufficiently small \( \varepsilon > 0 \), the scaled control \( u = \varepsilon \tilde{u} \) remains feasible with respect to \eqref{eq:QP_general_hdot} and \eqref{eq:QP_general_gdot}, achieves \( \dot{f}(z; u) < 0 \), and incurs a strictly lower cost than the trivial choice \( u = 0 \) or any other control $u$ yielding \( \dot{f}(z;u) = 0 \). This implies that the optimal control \( u^*_z \) must satisfy \( \dot{f}(z;u^*_z) < 0 \).

To prove the existence of such a control \( \tilde{u} \) at any non-stationary point \( z \), let \( J \) denote the set of active inequality constraints at \( z \), i.e.,
$
J = \{ j \in \{1, 2, \dots, k\} : h_j(z) = 0 \}.
$
Define the matrix \( A \in \mathbb{R}^{(1+|J|+2n) \times d} \) and vector \( \mathbf{b} \in \mathbb{R}^{1+|J|+2n} \) as:
\begin{align*}
A^\top &= \begin{bmatrix}
\nabla f(z) & -\nabla h_j(z)_{j \in J} & \nabla g_{i}(z)_{i = 1}^n & -\nabla g_{i}(z)_{i = 1}^n
\end{bmatrix},\\
\mathbf{b}^\top &= \begin{bmatrix}
-1 & 0 & \cdots & 0
\end{bmatrix}
\end{align*}
We now invoke \textit{Corollary} \ref{cor:Farkas}, a direct consequence of Farkas' Lemma \cite{10.5555/1214763}, both stated in Appendix~\ref{App:lemmas}, to show that no nonnegative vector $\mathbf{y} \in \mb R^{1+\abs{J}+2n}$ satisfies $A^\top \mathbf{y} = 0$ and $b^\top \mathbf{y} < 0.$ 
Suppose, for the sake of contradiction, \( \exists \mathbf{y} \) nonnegative, satisfying $A^\top \mathbf{y} = 0, \ b^\top \mathbf{y} < 0.$
Write the components of \( \mathbf{y} \) as \( \mathbf{y}^\top = \begin{bmatrix} y_0 & y_h^\top & y_{g^+}^\top & y_{g^-}^\top \end{bmatrix} \), where:
\begin{itemize}
    \item \( y_0 \in \mathbb{R}, \ y_h = [y_{h_{1}}, \dots, y_{h_{|J|}}]^\top, \in \mathbb{R}^{\abs{J}} \),
    \item \( y_{g^\pm} = [y_{g^\pm_1}, \dots, y_{g^\pm_n}]^\top \in \mb R^n \).
\end{itemize}
By expanding \(A^\top \mathbf{y} = 0\) and \(\mathbf{b}^\top \mathbf{y} < 0\), we get:
\begin{equation*} \label{eq: KKT}
\nabla f(z) y_0 - \sum_{j \in J} \nabla h_j(z) y_{h_j} + \sum_{i=1}^n \nabla g_i(z) (y_{g^+_i} - y_{g^-_i}) = 0,
\end{equation*}
and \(y_0 > 0\). Therefore, by dividing the above equation by \(y_0\), and denoting \(y_{h_j}/y_0\) by \(\mu_j \geq 0, \ \forall j \in J\), and \((y_{g^+_i} - y_{g^-_i})/y_0\) by \(\lambda_i\) \( \forall i = 1, \dots, n\), we get:
\[
\nabla f(z) - \sum_{j \in J} \nabla h_j(z) \mu_j + \sum_{i=1}^n \nabla g_i(z) \lambda_i = 0,
\]
which is precisely the KKT stationarity condition, contradicting the assumption that \(z\) is not a stationary point. Thus, the second statement in \textit{Corollary}~\ref{cor:Farkas} must hold—namely, there exists a solution \( \tilde{u} \in \mathbb{R}^b \) to the system \( A \tilde{u} \leq \mathbf{b} \). 
Consequently, \( \tilde{u} \) satisfies \eqref{eq:QP_general_hdot} for all \( j \in J \), as well as \eqref{eq:QP_general_gdot} and \( \dot{f}(z;\tilde{u}) \leq -1 < 0 \). 
For inactive constraints \( j \notin J \) (i.e., \( h_j(z) > 0 \)), we can scale \( \tilde{u} \) appropriately to ensure \eqref{eq:QP_general_hdot} holds for all \( j \). This completes the proof of part 1. \hfill \scalebox{0.9}{\( \blacksquare \)}

\subsubsection*{Part 2}

At each \( z \in \mathcal{C} \cap \mathcal{S} \), all conditions (1)–(5) of Theorem~1 in~\cite{6760327} are satisfied. Conditions (2)–(5) follow directly from the structure of the problem, while condition (1) holds due to the linear independence of the active constraint gradients, as stated in Assumption~\ref{assump:constraint_gradients}. Therefore, the solution map \( z \mapsto u^*(z) \) is locally Lipschitz continuous on \( \mathcal{C} \cap \mathcal{S} \). \hfill \scalebox{0.9}{\( \blacksquare \)}
\subsubsection*{Part 3}

Part 2 together with Assumption~\ref{assump:coercivity}, ensures existence and uniqueness of solutions to the closed-loop system \( \dot{z}(t) = u^*(z(t)) \) for all \( t \geq 0 \) and any initial condition \( z(0) \in \mathcal{C} \cap \mathcal{S} \)~\cite{khalil2002nonlinear}. As shown in~\cite{7782377}, constraint~\eqref{eq:QP_general_hdot} guarantees forward invariance of \( \mathcal{S} \), while~\eqref{eq:QP_general_gdot} ensures trajectories remain in \( \mathcal{C} \). Thus, \( \mathcal{C} \cap \mathcal{S} \) is forward-invariant under \(u^*(z)\). By LaSalle's Invariance Principle~\cite{khalil2002nonlinear}, trajectories converge to the largest invariant set where \( \dot{f}(z;u^*_z) = 0 \), which, as established in Part~1, coincides with the KKT points of \eqref{eq:opt_primal}.
\end{proof}

\subsection{SGF Aproach}
As discussed, SGF is an alternative way of viewing optimization problem \eqref{eq:opt_primal} as a control design task. 
In this approach, control inputs \( \rb{u,v} \in \mb R^k \times \mb R^n \) are designed for the following control-affine system:
\begin{align*}
    \dot z & = -\nabla f \rb z - \nabla h\rb{z}^\top u - \nabla g\rb z ^\top v, \quad z(0) \in \mathcal{C} \cap \mathcal{S} \\
    &\textrm{where }
        \begin{cases}
            h\rb{z} = \rcb{h_1\rb z \ \dots \ h_k\rb z}^\top, \\  
        g\rb{z} = \rcb{g_1\rb z \ \dots \ g_n\rb z}^\top.
        \end{cases}
\end{align*}
This system is interpreted as a gradient flow on \( f \), modified by control actions \( (u, v) \) that intervene to preserve feasibility when the state approaches constraint violation. The control inputs \( (u, v) \) are computed by solving the following QP:
\begin{align} \label{eq:QP_general_primal_dual}
    &\min_{\rb{u, v} \in \mb R^k_{\leq 0} \times \mb R^n} \big \|\nabla h\rb{z}^\top u + \nabla g\rb z ^\top v\big\|^2 \\ 
    & \quad \textrm{s.t.}
    \begin{cases}
        \nabla g \cdot \nabla g^\top u + \nabla g \cdot \nabla h^\top u 
        \leq -\nabla g \cdot \nabla f + \alpha g\rb z, \\ 
        \nabla h \cdot \nabla g^\top u + \nabla h \cdot \nabla h^\top u 
        = -\nabla h \cdot \nabla f + \alpha h\rb z,
    \end{cases} \notag
\end{align}
where $\alpha > 0$ is a constant.
The well-posedness of the SGF approach is established in \cite{sgf_allibhoy2023control}, where the flow is shown to be locally Lipschitz continuous, defined on an open set containing the feasible region, and converging to the KKT points of problem~\eqref{eq:opt_primal}.

\subsection{Instantiating our CBF-based Approach for \( \textrm{P}_{\textrm{constr}}(\beta) \)}


To ensure non-negativity of the objective \( \mathcal{F}^\beta \) in \( \textrm{P}_{\textrm{constr}}(\beta) \), as required by our CBF-based approach, we shift it by a constant and define the resulting function as
\begin{align}
\tilde{\mathcal{F}}^\beta := \frac{\log M}{\beta} + \sum_{i,j} p_i p_{j|i} \left( d(x_i, y_j) + \frac{\log p_{j|i}}{\beta} \right).  
\end{align}
Recall that \( \textrm{P}_{\textrm{constr}}(\beta) \) must be solved repeatedly at increasing values of $\beta_k$, with the solution of \( \textrm{P}_{\textrm{constr}}(\beta_k) \) serving as an initial condition for solving \( \textrm{P}_{\textrm{constr}}(\beta_{k+1}) \).
Consequently, for each \( \beta \) we introduce the following control system:
\begin{align*}
    &\dot{p}_{j|i}^{\beta} = v_{ij}, \,   && p_{j|i}^\beta(0) = p_{j|i}^0 \in (0,1)  \ &&& \forall i, j,\\
    &\dot{y}_j^\beta = u_j, \, && y_j^\beta \rb 0 = y_j^0,   &&& \forall j,
\end{align*}
Here, \( \{p_{j|i}^0\} \) and \( \{y_j^0\} \) represent a feasible initial condition. We first present the quadratic program used for control design, then verify that Assumptions~\ref{assump:constraint_gradients} and~\ref{assump:coercivity} hold for this problem. This enables the use of Theorem~\ref{Thrm:main_theorem} (part 3) to establish that \( \{p_{j|i}^\beta(t)\} \) and \( \{y_j^\beta(t)\} \) converge to a KKT point of \( \textrm{P}_{\textrm{constr}}(\beta) \) as \( t \to \infty \). For notational simplicity, we omit the superscript \( \beta \) in the remainder of this section.
\begin{align} \label{eq:QP_F}
    \min_{\substack{\{v_{ij}\}, \{u_j\},\\ (\delta, \varepsilon_k) \in \mathbb{R}^2}} 
    &\sum_{i,j} v_{ij}^2 + \sum_j \|u_j\|^2 + q_1 \, \delta^2 + q_2 \, \varepsilon_k^2 \\[2mm]
    \text{s.t.} \quad 
    &\dot{\tilde{\mathcal{F}}}(\{p_{j|i}\}, \{y_j\}) < - \mu \,\tilde{\mathcal{F}} + \delta, \notag\\[1mm]
    &\dot \phi(\{p_{j|i}\}_j) = 0, 
    &&\forall i, \notag\\[1mm]
    &\dot \psi_c \left(\{p_{j|i}\}_{i}\right) \geq -\alpha_{\psi_c} \psi_c \left(\{p_{j|i}\}_{i}\right), 
    &&\forall j, \notag \\[1mm]
    &\dot \psi_l \left(\{p_{j|i}\}_{i}\right) \geq -\alpha_{\psi_l} \psi_l \left(\{p_{j|i}\}_{i}\right), 
    &&\forall j, \notag \\[1mm]
    &\dot \xi(p_{j|i}) \geq -\alpha_{\xi} \xi(p_{j|i}), 
    &&\forall i,j. \notag
\end{align}
\noindent where 
\begin{list}{\arabic{enumi})}{\usecounter{enumi} \leftmargin=0em \labelsep=0.5em}
\item \( q_1, q_2, \mu, \alpha_{\psi_c}, \alpha_{\psi_l}, \) and \( \alpha_{\xi} \) are positive constants.
\item $\phi, \psi_c, \psi_l$ and $\xi$ are smooth functions as defined below:
\begin{subequations}        
    \begin{align}
        &\phi:\mathbb{R}^M \to \mathbb{R},  \ \phi(\{p_{j|i}\}_j) = \sum_j p_{j|i} - 1, \label{eq:phi}\\
        &\psi_c:\mathbb{R}^N \to \mathbb{R},  \ \psi_c\left(\{p_{j|i}\}_i\right) = C_j - \sum_i p_i p_{j|i} \, c_{ij}, \label{eq:psi_c}\\
        &\psi_l:\mathbb{R}^N \to \mathbb{R},  \ \psi_l\left(\{p_{j|i}\}_i\right) = \sum_i p_i p_{j|i} \, c_{ij} - L_j, \label{eq:psi_l}\\
        &\xi:\mathbb{R} \to \mathbb{R},  \ \xi(p_{j|i}) = p_{j|i}(1 - p_{j|i}). \label{eq:xi}
    \end{align}
\end{subequations}
Here, \eqref{eq:phi} denotes the equality constraint manifold, and \eqref{eq:psi_c}, \eqref{eq:psi_l}, and \eqref{eq:xi} act as CBFs. Any feasible \(\{p_{j|i}\}\) must satisfy: \(\phi(\{p_{j|i}\}_j) = 0, \forall \ i\), \(\psi_c(\{p_{j|i}\}_i) \geq 0\), \(\psi_l(\{p_{j|i}\}_i) \geq 0, \forall \ j\), and \(\xi(p_{j|i}) \geq 0, \forall \ i, j\). Time derivatives appearing in \eqref{eq:QP_F} are provided in Appendix~\ref{App:time_derivatives}. 
\end{list}

\begin{remark}
    Since \(\xi\) acts as a CBF, it ensures that \(p_{j|i}(t)\) remains strictly within \((0,1)\) \( \forall \ t \geq 0\), which in turn guarantees that \(\tilde{\mathcal{F}}\) remains $C^1$ over time. Lipschitz continuous partial derivatives are provided in Appendix~\ref{App:time_derivatives}.    
\end{remark}
\begin{remark}
    Enforcing minimum utilization constraints with \(L_j > 0\) for all \(j\) guarantees that if any resource location \(y_j\) tends to infinity (i.e., \(\|y_j\| \to \infty\)), then \(\tilde{\mathcal{F}} \to \infty\), thereby satisfying the coercivity assumption.
\end{remark}

\begin{remark}
    Linear independence among active constraint gradients fails if all resources have exactly one active capacity constraint. To prevent this, it suffices that one resource remains strictly within bounds. If for some \( k \in \{1, \dots, M\} \), the initial \( \{p_{k|i}^0\} \) satisfy 
    $
    \psi_c\big(\{p_{k|i}^0\}_i\big) > 0, \quad \psi_l\big(\{p_{k|i}^0\}_i\big) > 0,
    $
    then these remain strictly positive for all \( t \ge 0 \), so neither constraint for resource \( k \) becomes active, preserving linear independence of active constraints along the trajectory.
\end{remark}



\begin{remark}
At any KKT point of \( \textrm{P}_{\textrm{constr}}(\beta) \), the association probabilities \( \{p_{j|i}^\beta\} \) satisfy \eqref{cfpconst}, \eqref{relconstr}, while the facility locations \( \{y_j^\beta\} \), being unconstrained, must lie at the weighted centroids of their assigned demand points as given in \eqref{eq:unconstrained_updates}.
\end{remark}

\section{SIMULATIONS AND RESULTS} \label{sec:Simulations}

Our simulations demonstrate constraint handling, convergence rates, cost-effectiveness, and scalability on capacitated FLP using our approach alongside SGF \cite{sgf_allibhoy2023control}, SLSQP \cite{boggs1995sequential}, and the DA-based penalty (DA-P) method \cite{srivastava2022inequality}. We implement QPs for our method and SGF using the OSQP solver \cite{OSQP_stellato2020osqp} in \texttt{CVXPY}, implement corresponding state-space dynamics with an adaptive Euler-forward method \cite{malitsky2019adaptive}, and implement SLSQP via \texttt{SciPy} \cite{virtanen2020scipy}.
\vspace{-5 pt}
\subsection{Constraint handling, speed and cost:}
Consider a scenario with \( N = 400 \) users and \( M = 4 \) facilities, distributed in a $40\times 40$ rectangular area. The user clusters generated by normal distributions in the proportions given by $\cb{0.09, 0.32, 0.19, 0.40}$. The goal is to find an optimal pair of facility locations and facility-user assignments such that the total servicing cost is minimized. Additionally, the facilities have upper capacity constraints given by $\cb{0.20, 0.41, 0.27, 0.20}$, denoting the fraction of users that each facility can service. The square Euclidean distance serves as a servicing cost measure for each user-facility pair. 
\begin{figure}[tbhp]
    \centering
    \subfigure[Our Approach]{\includegraphics[trim={1.6cm 0.9cm 2.7cm 2.2cm},clip,width=0.45\columnwidth]{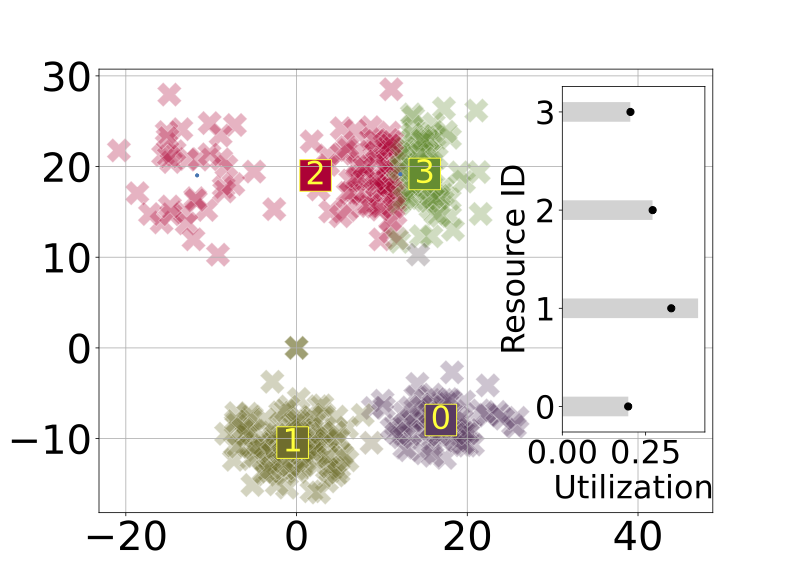}}
    \subfigure[SGF]{\includegraphics[trim={1.6cm 0.9cm 2.7cm 2.2cm},clip,width=0.45\columnwidth]{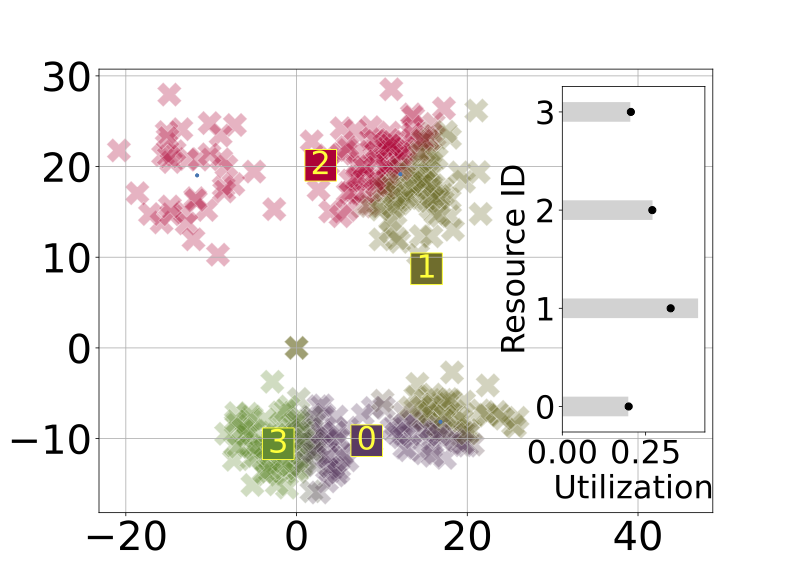}}    
    \subfigure[SLSQP]{\includegraphics[trim={1.6cm 0.9cm 2.7cm 2.2cm},clip,width=0.45\columnwidth]{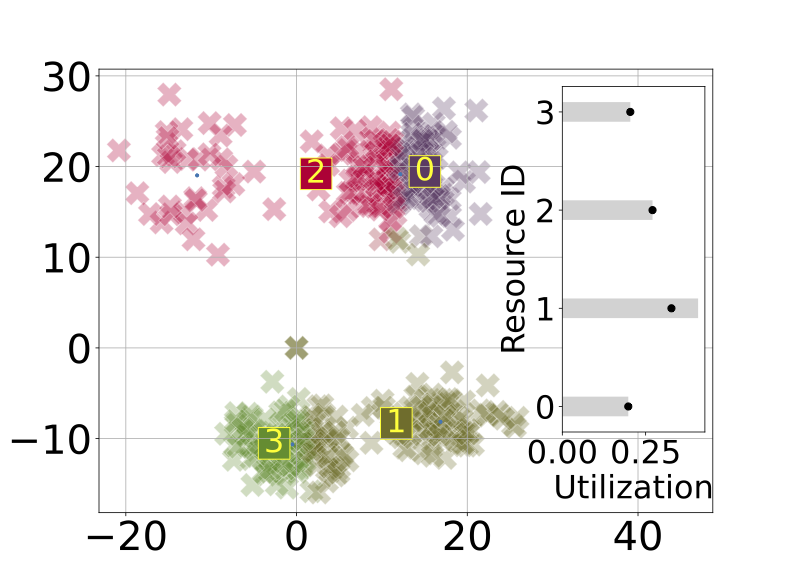}} 
    \subfigure[DA-P]{\includegraphics[trim={1.6cm 0.9cm 2.7cm 2.2cm},clip,width=0.45\columnwidth]{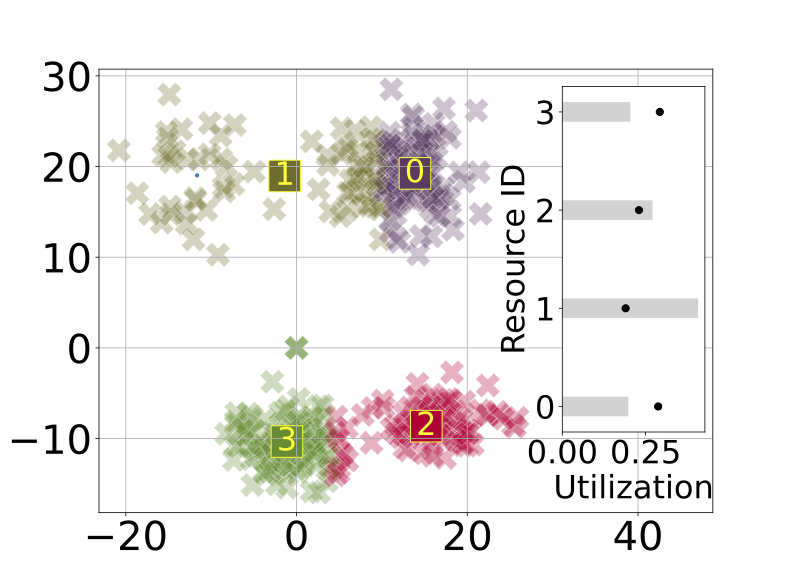}}
    \subfigure[Runtime and cost comparison as $\beta$ grows $10^{-3} \rightarrow 10^2$.]{\includegraphics[trim={1.1cm 0.08cm 3.0cm 1.8cm},clip,width=0.9\columnwidth]{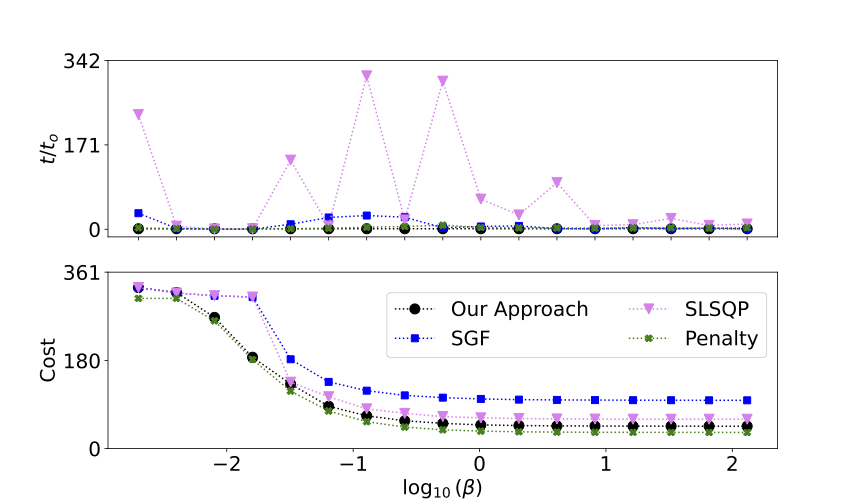}}
    \vspace{-5 pt}
    \caption{\small The figure shows a capacitated FLP with 400 demand points in 4 clusters, solved using the four methods. Final resource utilization is shown to the right of each subplot. All the approach maintain feasibility except the DA-based penalty method. Runtimes (in sec): \(\{46, 210, 1600, 60\}\), Final costs: \(\{46, 99, 60, 33\}\) units.}
    \label{fig:compare_SGF_Our_AmberPenalty}
\end{figure}
Figure~\eqref{fig:compare_SGF_Our_AmberPenalty} shows the solution obtained by aforementioned approaches. Each facility and its assigned users share the same color. The final capacity of each facility is shown at the right of each subplot. Figure~(\ref{fig:compare_SGF_Our_AmberPenalty}e) compares the runtime and corresponding cost as $\beta$ is increased from $10^{-3}$ to $100$ during annealing iterations. We observe, the solution speed of our approach to perform the entire annealing process is comparable to DA-P, and $5$ and $35$ times faster than SGF and SLSQP respectively. Further, our approach results $\sim 25 \%$ and $\sim 53 \%$ lesser costs w.r.t. SLSQP and SGF respectively. Although, DA-P yields the lowest cost, it violates the utilization constraints. Further, DA-P demands careful, problem-specific hyperparameter tuning to prevent cost distortion or numerical instability. 

Next, we numerically compare convergence time and the final cost for solving  \( \textrm{P}_{\textrm{constr}}(\beta) \) using our approach, SGF and SLSQP with growing problem size.
\subsection{Convergence rates at fixed $\beta$:}
For a 2D FLP with $N$ demand points and $M$ facilities, we define problem size for \( \textrm{P}_{\textrm{constr}}(\beta) \) by $NM + 2M,$ representing the total number of decision variables. Table~\ref{tab:fix_beta_compare} reports the corresponding runtime and cost at $\beta = 0.001, 1$ and $100$, as problem size grows from $204$ to $3010$. For each $\beta$ and $\rb{N,M}$, an FLP instance is generated by selecting user locations, initial facility locations and facility utilization limits at random. Further, \( \textrm{P}_{\textrm{constr}}(\beta) \) is solved by starting at same initial conditions for the three methods. 

We observe that our approach and SLSQP yield identical cost values across all instances, where as SGF results higher costs particularly for larger problem instances at $\beta=1$ and $100$. Furthermore, as problem size increases from $204$ to $3010,$ the runtime of our approach grows by a factor of $\sim16$, compared to $\sim110$ for SGF and an astonishing $\sim13\times 10^3$ for SLSQP. An average across all the simulation scenarios in Table~\ref{tab:fix_beta_compare} shows that our approach is $20\times$ and $240\times$ faster as compared to SGF and SLSQP respectively. These observations show that our approach remains cost effective with minimal scaling as compared to SGF and SLSQP. 
In the next subsection, we demonstrate the capability of our approach in handling a significantly large FLP problem.
\begin{table}[tbhp]
    \centering
    \begin{tabular}{|c|c|cc|cc|cc|}
        \hline
        \multirow{2}{*}{$\beta$} & \multirow{2}{*}{P-Size} & \multicolumn{2}{c|}{Our Method} & \multicolumn{2}{c|}{SGF} & \multicolumn{2}{c|}{SLSQP} \\
        \cline{3-8}
        &  & Cost & Time & Cost & Time & Cost & Time \\
        \hline
        
        \multirow{3}{*}{$10^{-3}$} 
        & 204   & 37.6  & 0.3 & 36.3 & 0.8 & 38.1 & 0.5 \\
        & 404   & 46.2  & 0.4 & 44.9 & 1.2 & 45.9 & 2.0 \\
        & 906   & 51.8  & 0.7 & 52.0 & 11.3 & 51.9 & 35.8 \\
        & 1608  & 63.2  & 1.5 & 63.2 & 43.2 & 63.2 & 271.5 \\
        & 2008  & 52.8  & 1.2 & 52.6 & 62.1 & 52.6 & 489.2 \\
        & 3010  & 118.5 & 3.4 & 118.4 & 163.8 & 118.3 & 2205.5 \\
        \hline
        
        \multirow{3}{*}{1} 
        & 204  & 10.8 & 0.4 & 11.7 & 2.2 & 10.8 & 0.3 \\
        & 404  & 12.0 & 0.6 & 16.5 & 4.3 & 12.1 & 2.9 \\
        & 906  & 12.5 & 1.1 & 23.3 & 14.8 & 12.3 & 66.0 \\
        & 1608  & 10.3 & 1.5 & 22.9 & 51.8 & 10.5 & 356.7 \\
        & 2008  & 10.8 & 3.7 & 26.5 & 81.8 & 12.2 & 959.8 \\
        & 3010  & 11.5 & 9.9 & 24.0 & 170.8 & 11.0 & 5670.1 \\
        \hline
        
        \multirow{3}{*}{100} 
        & 204  & 10.3 & 0.3 & 10.9 & 3.8 & 10.3 & 0.3 \\
        & 404  & 12.0 & 0.5 & 12.9 & 4.1 & 12.0 & 2.3 \\
        & 906  & 10.9 & 1.0 & 21.6 & 15.7 & 10.9 & 11.0 \\
        & 1608  & 12.0 & 1.3 & 20.6 & 50.3 & 11.5 & 364.0 \\
        & 2008  & 15.0 & 1.6 & 22.5 & 81.7 & 14.7 & 765.4 \\
        & 3010  & 13.7 & 4.2 & 27.5 & 185.1 & 13.0 & 4570.7 \\
        \hline
    \end{tabular}
    \caption{\small Time (sec) \& cost comparison for solving \( \textrm{P}_{\textrm{constr}}(\beta) \) as problem size (P-size) grows.}
    \label{tab:fix_beta_compare}
\end{table}
\vspace{-20 pt}
\subsection{Solving a large scale problem:}
\begin{figure}[tbhp!]
    \centering
    \includegraphics[trim={6.0cm 4.8cm 8.0cm 8.0cm}, clip, width=0.9\columnwidth]{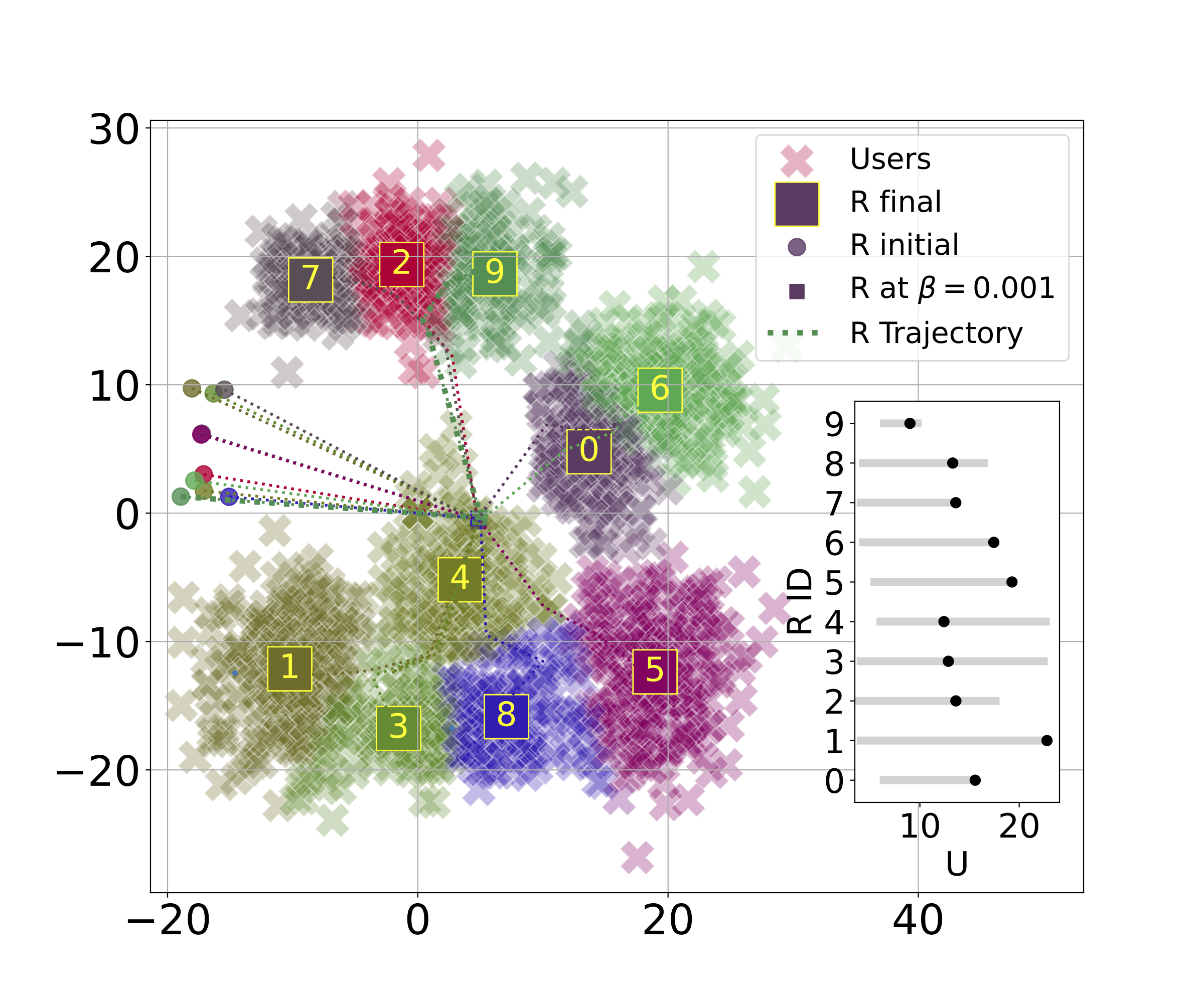}
    \vspace{-6 pt}
    \caption{\small Capacitated FLP solution using our CBF-based approach for $N=1000, M=10.$ The cluster split of users: $\rcb{0.11, 0.07, 0.11, 0.09, 0.13, 0.14, 0.02, 0.11, 0.14, 0.08}$ and facility (R) utilization (U) constraints are shown at the bottom right. The figure also shows splitting of facilities into distinct clusters as $\beta \in \rcb{10^{-3},100}$ is increased during annealing.}
    \label{fig:cbf_2000_10}
    \includegraphics[trim={2.4cm 16.0cm 7.0cm 2.0cm},clip,width=0.9\columnwidth]{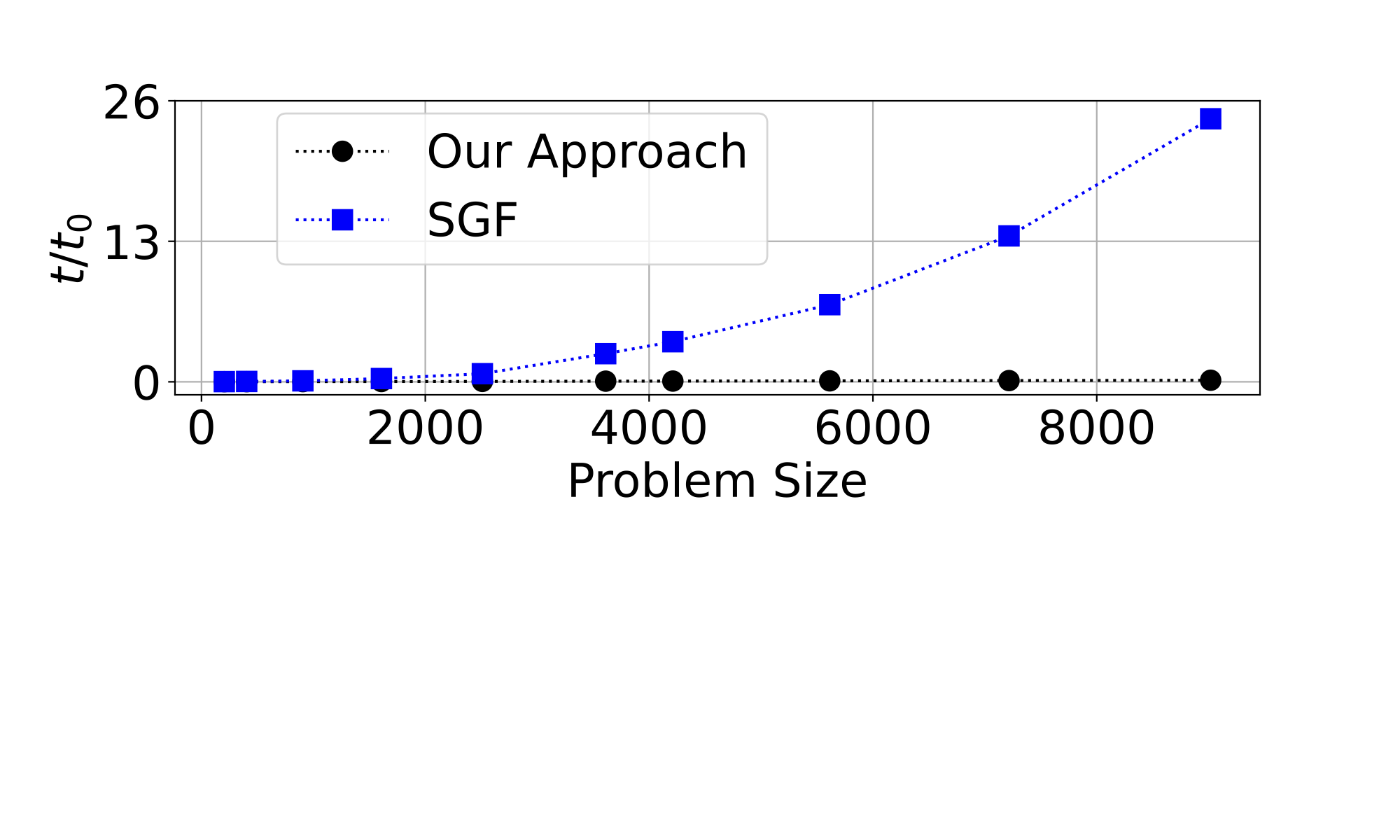}
    \vspace{-8 pt}
    \caption{\small Runtime comparison for solving QP in our approach \eqref{eq:QP_general} and SGF \eqref{eq:QP_general_primal_dual} as problem size grows for fixed $\beta$ and state. The SGF QP runtime ($t$) is shown relative to our QP runtime ($t_0$). Runtime is reported as an average for 10 random $\beta$ instances, while using the same state for both methods at each $\beta$.}
    \label{fig:Our_SGF_control_compare_2025_3_29_19_54_51}
\end{figure}

We consider a large scale 2D FLP consisting of $N=2000$ users to be serviced by $M=10$ facilities, with upper and lower utilization constraints on each facility as shown in Figure~\eqref{fig:cbf_2000_10}. The decision space for this problem belongs to $\mb R^{20} \times \cb{0,1}^{20000}.$ Our algorithm solves the problem in $19$ minutes and the final cost is $18.1$ units. We also attempted to solve the problem using SGF, which required more than 60 minutes to minimize free energy for each $\beta.$ 
Although both the methods are CBF-CLF based, the difference in convergence rates to achieve KKT conditions is attributed to the rate at which control actions are computed via QPs, where the SGF QP \eqref{eq:QP_general_primal_dual} scales poorly as compared to QP for our approach \eqref{eq:QP_F} as observed in Figure~\eqref{fig:Our_SGF_control_compare_2025_3_29_19_54_51}.

These results demonstrate the proficiency of our algorithm in constraint handling constraints, scalability with respect to time and costs for large-scale resource allocation problems. 
\section{Conclusion and Future works}
We addressed capacity-constrained resource allocation problems involving joint optimization of resource attributes (e.g., spatial positions) and demand assignments. The problem is NP-hard and non-convex, and inequality constraints further increase its complexity. Leveraging the MEP-based DA framework, we reformulated the inner-loop optimization of free energy as a control design problem, using CLFs and CBFs to ensure descent toward stationarity and constraint satisfaction. We proved convergence to a KKT point under mild assumptions on the initial conditions. Simulations show that our CBF-based method enforces constraints, yields competitive costs, shows negligible growth in convergence times compared to benchmarks, and scales to problem sizes previously intractable with standard solvers.

As future work, we aim to extend our approach to dynamic settings where demand points evolve over time according to prescribed dynamics \( \dot{x}_i = \kappa(\{x_i\}) \), requiring the resource controllers to track moving local minima. Another promising direction is a theoretical investigation into the quality of the obtained solutions, including strategies such as structured perturbations to escape poor local minima.

\section{APPENDIX}

\subsection{Supporting Lemmas} \label{App:lemmas}

\begin{lemma}[Farkas' Lemma, Variant (iii)]
Let \( A \in \mathbb{R}^{s \times r} \) be a real matrix, and \( \mathbf{b} \in \mathbb{R}^s \). The system \( A \mathbf{x} \leq \mathbf{b} \) has a solution if and only if every nonnegative \( \mathbf{y} \in \mathbb{R}^s \) such that \( A^\top \mathbf{y} = 0 \) also satisfies \( \mathbf{b}^\top \mathbf{y} \geq 0 \). \textit{Proof:} section 6.4 of \cite{10.5555/1214763}.
\end{lemma}

\begin{corollary} \label{cor:Farkas}
Exactly one of the following must be true: 
\begin{enumerate}
    \item There exists a nonnegative \( \mathbf{y} \in \mathbb{R}^s \) such that \( A^\top \mathbf{y} = 0 \) and \( \mathbf{b}^\top \mathbf{y} < 0 \).
    \item The system \( A \mathbf{x} \leq \mathbf{b} \) has a solution for \( \mathbf{x} \in \mathbb{R}^r \).
\end{enumerate}
\end{corollary}

\subsection{Time Derivatives of the CLF and CBFs} \label{App:time_derivatives}
\vspace{-10pt}
{\small
\begin{align*}
    &\dot{\tilde{\mathcal{F}}}(\{p_{j|i}\}, \{y_j\}) = \sum_{i,j} \nabla \tilde{\mathcal{F}}_{p_{j|i}} v_{ij} + \sum_j \left \langle \nabla \tilde{\mathcal{F}}_{y_j}, u_j \right \rangle,\\
    & \nabla_{y_j} \tilde{\mathcal{F}} = 2 \sum_i p_i p_{j|i} (y_j - x_i), \quad  \forall j,\\
    & \nabla_{p_{j|i}} \tilde{\mathcal{F}} = p_i \left( \| x_i - y_j \|^2 + \frac{1}{\beta} ( \log p_{j|i} + 1 ) \right)  \quad \forall i,j. \\
    &\dot \psi \left(\{p_{j|i}\}_{i}\right) = \pm \sum_i p_i v_{ij} c_{ij}, \text{with } + \text{ for } \psi_l, \text{ and } - \text{ for } \psi_c,\\
    &\dot{\phi} (\{p_{j|i}\}_j) = \sum_j v_{ij}, \quad \text{and} \quad \dot \xi(p_{j|i}) = (1 - 2p_{j|i}) v_{ij}.
\end{align*}
}

\bibliographystyle{IEEEtran}
\bibliography{mybibfile}

\end{document}